\documentclass[11pt]{article}
\usepackage{geometry}                
\usepackage{wrapfig}
\usepackage{graphicx}
\usepackage{amssymb}
\usepackage{amsmath}
\usepackage{url}

\DeclareMathOperator{\sing}{sing}
\DeclareMathOperator{\cosg}{cosg}
\DeclareMathOperator{\tang}{tang}

\usepackage{epstopdf}
\usepackage[toc,page]{appendix}
\title{Non-Euclidean Triangle Centers}
\author{Robert A. Russell}

\begin{document}
\def\sing{\mathop{\rm sing}\nolimits}
\def\cosg{\mathop{\rm cosg}\nolimits}
\def\tang{\mathop{\rm tang}\nolimits}
\maketitle

\begin{abstract}
Non-Euclidean triangle centers can be described using homogeneous coordinates that are proportional to the generalized sines of the directed distances of a given center from the edges of the reference triangle.  Identical homogeneous coordinates of a specific triangle center may be used for all spaces of uniform Gaussian curvature.  We also define the {\em median point} for a set of points in non-Euclidean space and a {\em planar center of rotation} for a set of points in a non-Euclidean plane.
\end{abstract}

\section{Introduction}

Clark Kimberling's on-line {\em Encyclopedia of Triangle Centers}~\cite{Kimberling} is a collection of thousands of Euclidean triangle centers.  It provides descriptions and trilinear coordinates for each center, along with additional information.  There does not appear to be a similar collection for non-Euclidean triangle centers, which can also be given similar coordinate ratios.

In {\em Non-Euclidean Geometry}~\cite{Coxeter}, H. S. M. Coxeter describes the use of homogeneous coordinates for non-Euclidean spaces of uniform Gaussian curvature.  Coxeter mentions the homogeneous coordinates of three triangle centers.  These are the incenter ($1:1:1$), the orthocenter ($\sec A:\sec B:\sec C$), and the intersection of the medians ($\csc A:\csc B:\csc C$).  He notes that these are the same as the Euclidean trilinear coordinates.

\section{Euclidean coincidence}

In Euclidean geometry there is more coincidence for triangle centers than in non-Euclidean geometry.  For example, in Euclidean geometry, the center of rotation of the vertices, the center of rotation of the triangle, the intersection of the medians, and the intersection of the area-bisecting cevians are the same point.  In non-Euclidean geometry these are all distinct.  Although their trilinear coordinates in Euclidean geometry are identical, in non-Euclidean geometry the homogeneous coordinates are not.

In non-Euclidean space, the homogeneous coordinates of the circumcenter, with $S=(A+B+C)/2$, are ($\sin(S-A):\sin(S-B):\sin(S-C)$).  In Euclidean space, $S$ must be $\pi/2$, so that the trilinear coordinates for the circumcenter are also ($\cos A:\cos B:\cos C$), since $\sin(\pi/2-A) = \cos A$.  But in non-Euclidean geometry, the triangle center with homogeneous coordinates ($\cos A:\cos B:\cos C$) is not the circumcenter of the triangle.  It is a different center, which happens to coincide with the circumcenter in Euclidean space.

\section{Generalized trigonometric functions}

We shall employ generalized trigonometric functions, since they apply to all spaces of uniform Gaussian curvature.
For a space of uniform Gaussian curvature $K$, the generalized sine function is defined as
\begin{eqnarray*}
	\sing(x) & = & \sum_{i=0}^\infty (-K)^i\frac{x^{2i+1}}{(2i+1)!}
			= {\sin x\sqrt{K}\over\sqrt{K}}
			= {\sinh x\sqrt{-K}\over\sqrt{-K}} \\
		& = & x -Kx^3/3!+K^2x^5/5!-K^3x^7/7!+\dots \\
		& = & x \hspace{10pt}\mbox{if $K=0$} \\
		& = & \sin(x) \hspace{10pt}\mbox{if $K=1$} \\
		& = & \sinh(x) \hspace{10pt}\mbox{if $K=-1$}.
\end{eqnarray*}
We pronounce $\sing(x)$ the same as ``singe x.''  It allows us to express the law of sines for any space of uniform Gaussian curvature $K$ as
\[
	\frac{\sing a}{\sin A} = \frac{\sing b}{\sin B} = \frac{\sing c}{\sin C}.
\]

The generalized cosine function is defined as
\begin{eqnarray*}
	\cosg(x) & = & \sum_{i\geq0}{(-K)^i\over(2i)!}x^{2i}
			= \cos x\sqrt{K}
			= \cosh x\sqrt{-K} \\
		& = & 1 -Kx^2/2!+K^2x^4/4!-K^3x^6/6!+\dots \\
		& = & 1 \hspace{10pt}\mbox{if $K=0$} \\
		& = & \cos(x) \hspace{10pt}\mbox{if $K=1$} \\
		& = & \cosh(x) \hspace{10pt}\mbox{if $K=-1$}.
\end{eqnarray*}
We pronounce $\cosg(x)$ as if it rhymed with ``dosage x.''  We can show that
\[
	\cosg^2x+K\sing^2x=1.
\]

Finally, the generalized tangent function is defined as
\begin{eqnarray*}
	\tang(x) & = & {\sing x\over\cosg x}
			= {\tan x\sqrt{K}\over\sqrt{K}}
			= {\tanh x\sqrt{-K}\over\sqrt{-K}} \\
		& = & x \hspace{10pt}\mbox{if $K=0$} \\
		& = & \tan(x) \hspace{10pt}\mbox{if $K=1$} \\
		& = & \tanh(x) \hspace{10pt}\mbox{if $K=-1$}.
\end{eqnarray*}
We pronounce $\tang(x)$ as if it rhymed with ``flange x.''

\section{Homogeneous coordinates}

\begin{figure}
	\begin{center}
	\begin{picture}(400,120)(20,0)
		\includegraphics[width=6in]{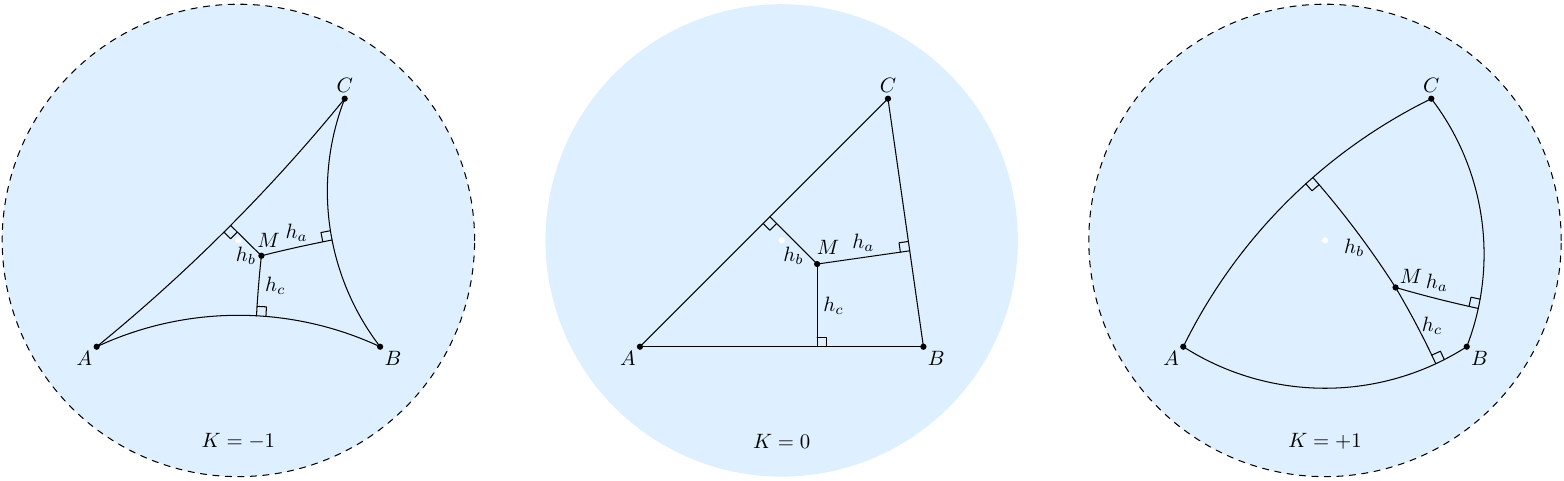}
	\end{picture}
	\end{center}
	\caption{Regardless of the Gaussian curvature $K$, the homogeneous coordinates of the point $M$ are $(\sing h_a:\sing h_b:\sing h_c)$, where $h_a$, $h_b$, and $h_c$ are the directed distances from $M$ to the corresponding sides of the reference triangle $ABC$.}
	\label{fig:centers1}
\end{figure}

Coxeter describes the homogeneous coordinates of a point $x$ as a triple ratio $(x_0:x_1:x_2)$.  These coordinates are equivalent if multiplied by the same value, so that $(x_0:x_1:x_2)=(\lambda x_0:\lambda x_1:\lambda x_2)$ for any nonzero number $\lambda$, in the same way that we can multiply the numerator and denominator of a fraction by the same value to obtain an equivalent fraction.  Given a reference triangle $ABC$ in a space of uniform Gaussian curvature $K$, we can obtain the homogeneous coordinates of a point by first calculating the directed distances of the point from the edges of the reference triangle.  The distances are positive (negative) if the point is on the same (opposite) side of the triangle edge as the remaining vertex.  As shown in Fig.~\ref{fig:centers1}, the directed distances of the point $M$ from the edges $a$, $b$, and $c$ of the reference triangle $ABC$ are $h_a$, $h_b$, and $h_c$ respectively.  The homogeneous coordinates of $M$ are then $(\sing h_a:\sing h_b:\sing h_c)$.  Since $\sing h = h$ in Euclidean geometry, these homogeneous coordinates are equivalent to trilinear coordinates when the Gaussian curvature $K=0$.

As with trilinear coordinates, we can also specify the homogeneous coordinates of a line $Y$ as a triple ratio $[Y_0:Y_1:Y_2]$.  A point $x$ is on a line $Y$ only if $\{xY\}=x_0Y_0+x_1Y_1+x_2Y_2=0$, so that the line through points $p$ and $q$ is $[p_1q_2-p_2q_1,p_2q_0-p_0q_2,p_0q_1-p_1q_0]$ as shown in $\mathsection4.3$ in Coxeter.  The homogeneous coordinates of the vertices of the reference triangle are $A=(1:0:0)$, $B=(0:1:0)$, and $C=(0:0:1)$.  The sides are $a=[1:0:0]$, $b=[0:1:0]$, and $c=[0:0:1]$.

Let us prove that these homogeneous coordinates are those defined by Coxeter.  In equation 12.14, Coxeter states that the distance from a point $(x)$ to a line $[Y]$ is 
\[
	\sin^{-1} \frac{|\{xY\}|}{\sqrt{(xx)}\sqrt{[YY]}}\hspace{20pt}\mbox{for $K=1$}\hspace{40pt}
	\sinh^{-1} \frac{|\{xY\}|}{\sqrt{(xx)}\sqrt{-[YY]}}\hspace{20pt}\mbox{for $K=-1$}
\]
Then the distances from a point $x$ to the three edges of the reference triangle are
\[
	\sing^{-1}\frac{x_0}{\sqrt{(xx)}\sqrt{K[YY]}}\hspace{25pt}
	\sing^{-1}\frac{x_1}{\sqrt{(xx)}\sqrt{K[YY]}}\hspace{25pt}
	\sing^{-1}\frac{x_2}{\sqrt{(xx)}\sqrt{K[YY]}}
\]
Taking the generalized sine of these distances gives us the homogeneous coordinates, since the denominators are equivalent.  Thus, the generalized sines of the directed distances of a point from the edges of the reference triangle are identical to the homogeneous coordinates described by Coxeter.

These homogeneous coordinates for points and lines allow us to do linear operations in spaces of constant curvature.  As described above, we can tell if a point is on a line and determine a line passing through two points.  Three points ($x$, $y$, and $z$) are collinear if
\[
	\begin{vmatrix}
		x_0 & x_1 & x_2 \\
		y_0 & y_1 & y_2 \\
		z_0 & z_1 & z_2 \\
	\end{vmatrix} = 0.
\]
The intersection point of two lines $X$ and $Y$ is
\[
	\left(X_1Y_2-X_2Y_1:X_2Y_0-X_0Y_2:X_0Y_1-X_1Y_0\right).
\]

\section{Finding the homogeneous coordinates of a triangle center}

\begin{wrapfigure}{r}{6.15cm}\includegraphics[scale=0.75]{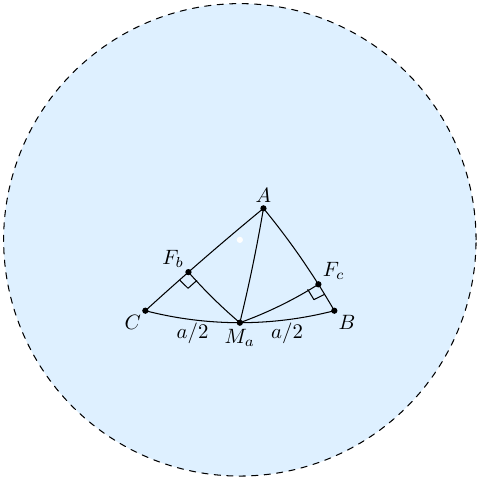} 
	\caption{What are the homogeneous coordinates of midpoint $M_a$?}
	\label{fig:cir_curv}
 \end{wrapfigure}

Let us find the homogeneous coordinates of the intersection of the medians of a non-Euclidean triangle.  First we want to find the homogeneous coordinates of the midpoint of an edge of the triangle.  Let $M_a$ be the midpoint of edge $a$.  Let $F_b$ be on edge $b$ such that $M_aF_b$ is orthogonal to edge $b$, and let $F_c$ be on edge $c$ such that $M_aF_c$ is orthogonal to edge $c$.  Then the homogeneous coordinates of $M_a$ are $(0:\sing M_aF_b:\sing M_aF_c)$.  Using the generalized law of sines on right triangle $M_aF_bC$, we have
\[
	{\sing M_aF_b\over \sin C}={\sing {a\over2}\over 1}.
\] 
Similarly,
\[
	{\sing M_aF_c\over \sin B}={\sing {a\over2}\over 1}.
\]
Substituting, we have the homogeneous coordinates
\[
	M_a=(0:\sing M_aF_b:\sing M_aF_c)=(0:\sin C\sing{a\over2}:\sin B\sing{a\over2})
	=(0:\sin C:\sin B).
\]
Similarly, $M_b=(\sin C:0:\sin A)$ and $M_c=(\sin B:\sin A:0)$.  Now we can determine the homogeneous coordinates of the medians $AM_a$, $BM_b$, and $CM_c$ to be
\[
	AM_a=[0:-\sin B:\sin C]\hspace{20pt}BM_b=[\sin A:0:-\sin C]\hspace{20pt}CM_c=[-\sin A:\sin B:0]
\]
The point $M$ at the intersection of any two of these medians is
\[
	M=(\sin B\sin C:\sin A\sin C:\sin A\sin B)=(\csc A:\csc B:\csc C).
\]
The use of the linear algebra of homogeneous coordinates makes this exercise easy.

\section{Coordinate conversion}

The stereoscopic projection of a space of uniform curvature projects points on an embedded curved space onto a Euclidean subspace, which we shall make tangent to the embedded curved space.  We place the origin at the tangent point of the curved space, and we know the Gaussian curvature $K$ of the curved space.  We shall want to convert coordinates of points on the projection to and from points on the embedded curved space.  We shall also want to convert homogeneous coordinates to and from these coordinates as well.

Consider a space of uniform curvature of dimension $d$.  We can assign a point on its stereographic projection ordinary Cartesian coordinates such as $(x_1,x_2,\dots,x_d)$ or polar coordinates such as $r\theta$, where $r$ is the distance of the point from the origin and $\theta=(\theta_1,\theta_2,\dots,\theta_d)$ is a sequence of direction cosines such that $1=\sum_i\theta_i^2$ and $r\theta_i=x_i$.  For points on the embedded curved space, we require an additional coordinate $x_0$, which will have an imaginary value in the case where $K<0$.

\begin{figure}
	\begin{center}
	\begin{picture}(400,200)(-36,0)
		\includegraphics[width=5in]{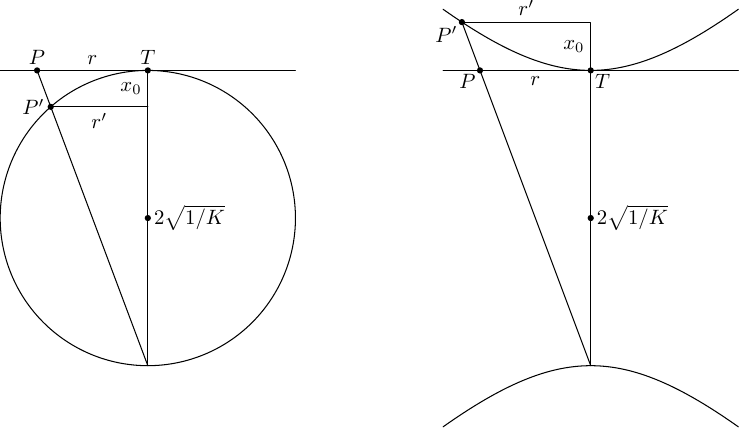}
	\end{picture}
	\end{center}
	\caption{Transforming the coordinates of points $P$ and $P'$ with positive (left) and negative (right) curvature.  The coordinate origin point $T$ is the tangent point of the curved surface.}
	\label{fig:centers2}
\end{figure}

In Fig.~\ref{fig:centers2}, we see cross sections of the embedded space and the image space on which it is projected, a horizontal line.  The coordinates of $P$, the stereographic projection, are $r\theta$.  We can convert these to the coordinates of $P'$, which are $(x_0,r'\theta)$.  The equations are
\[
	{r\over2\sqrt{1/K}}={r'\over2\sqrt{1/K}+x_0}\hspace{30pt}r'^2+(\sqrt{1/K}+x_0)^2=1/K
\]
Solving, we have
\[
	x_0={-2r^2\over(4+Kr^2)\sqrt{1/K}}\hspace{20pt}r'={4r\over4+Kr^2}
	\hspace{20pt}r=2\sqrt{-x_0\over Kx_0+2\sqrt{K}}\hspace{20pt}r={2-2\sqrt{1-Kr'^2}\over Kr'}
\]
Note that $r'$ is less than $r$ when $K<0$ and greater than $r$ when $K<0$.  Also the first coordinate $x_0$ is negative when $K>0$ and positive imaginary when $K<0$.

Conversion of the homogeneous coordinates is more complicated.  We first need the coordinates of the vertices of the triangle of reference in the embedded space.  We then find the equations of the three planes that pass through the edges of the reference triangle and the center of the sphere.  Next we translate each plane in the direction of its normal (positive is in the direction of the remaining vertex of the reference triangle) by the amount of the corresponding homogeneous coordinate.  These planes will intersect at a point, which we can project onto the sphere from its center.

\begin{figure}
	\begin{center}
	\begin{picture}(400,200)(0,0)
		\includegraphics[width=5in]{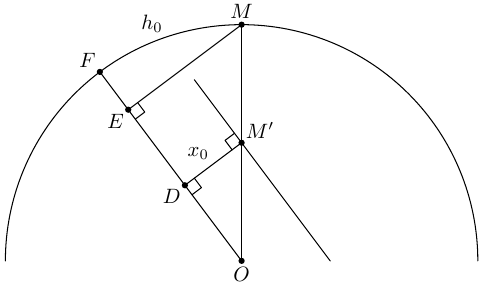}
	\end{picture}
	\end{center}
	\caption{A cross section through its center $O$ of an embedded sphere with radius $\sqrt{1/K}$.  Point $M$ is a triangle center, and point $F$ is the foot of a perpendicular drawn on the sphere from the center to an edge.}
	\label{fig:centers7}
\end{figure}

Why does this method work?  In Fig.~\ref{fig:centers7}, we have drawn a cross section of an embedded sphere with radius $\sqrt{1/K}$.  Point $F$ is the foot of a perpendicular drawn from a triangle center $M$ to an edge of the triangle.  Letting the homogeneous coordinates of $M$ be $(x_0:x_1:x_2)$, we draw a plane parallel to the plane through the edge at $F$ and the center $O$ such that the translation $DM'$ is equal to $x_0$.  Point $M'$ is the intersection of the plane with line $OM$.  The length of the arc $MF$ is $h_0$, and we can see that
\[
	\sin{h_0\over OM}={EM\over OM},
\]
which means that $EM=\sing h_0$ by definition.  From similar triangles, we can see that
\[
	{EM\over DM'}={\sing h_0\over x_0}={OM\over OM'}.
\]
But this will be true for the other edges as well, so that $(x_0:x_1:x_2)=(\sing h_0:\sing h_1:\sing h_2)$.  This means that the intersection of the planes at $M'$ can be centrally projected to obtain the center $M$.  Let's do examples with positive and negative Gaussian curvature.

\subsection{Converting homogeneous coordinates when $\mathbf{K=+1}$.}
\begin{figure}
	\begin{center}
	\begin{picture}(400,150)(0,0)
		\includegraphics[width=5in]{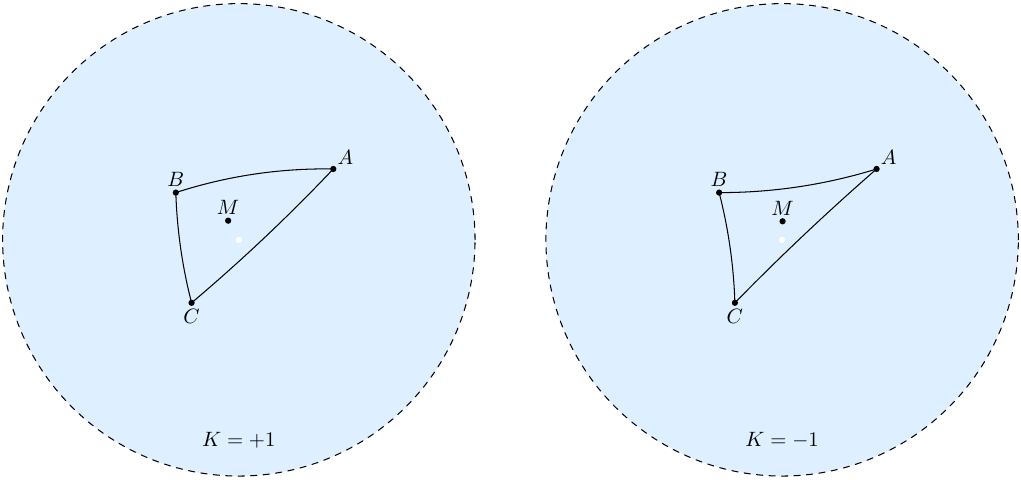}
	\end{picture}
	\end{center}
	\caption{The point $M$ is the median point of triangle $ABC$.}
	\label{fig:centers5}
\end{figure}

Let $K=1$ and the projected vertices of the triangle of reference be $A(4/5,3/5)$, $B(-8/15,2/5)$, and $C(-2/5,-8/15)$.  Then the polar coordinates are $r_A=1$, $\theta_A=(4/5,3/5)$, $r_B=2/3$, $\theta_B=(-4/5,3/5)$, $r_C=2/3$, and $\theta_C=(-3/5,-4/5)$.  Then we can calculate the coordinates on the embedded sphere using our formulas.  We have
\[
	A'=\left(-{2\over5},{16\over25},{12\over25}\right)
	\hspace{30pt}B'=\left(-{1\over5},-{12\over25},{9\over25}\right)
	\hspace{30pt}C'=\left(-{1\over5},-{9\over25},-{12\over25}\right).
\]
Since the center $O$ of the embedded sphere is at $(-\sqrt{1/K},0,0)=(-1,0,0)$, the three planes in which the triangle edges are embedded are
\begin{eqnarray*}
	OB'C': & & 15x_0+28x_1+4x_2+15=0 \\
	OC'A': & & 12x_0-60x_1+65x_2+12 = 0 \\
	OA'B': & & 288x_0+105x_1-500x_2+288=0
\end{eqnarray*}
The normals to these planes allow us to compute the cosines of the angles of the triangle of reference, which are
\[
	\cos A={35344\over29\sqrt{3259321}}
	\hspace{36pt}\cos B=-{1052\over29\sqrt{16769}}
	\hspace{36pt}\cos C={248\over\sqrt{326729}}
\]
From these we can compute the sines of the angles.
\[
	\sin A={38625\over29\sqrt{3259321}}
	\hspace{36pt}\sin B={3605\over29\sqrt{16769}}
	\hspace{36pt}\sin C={515\over\sqrt{326729}}
\]
We want to translate the planes above by the homogeneous coordinates, for which we shall use the intersection of the medians, ($\csc A:\csc B:\csc C$).  The new planes are
\begin{eqnarray*}
	15x_0+28x_1+4x_2 & = & \sqrt{15^2+28^2+4^2}\csc A-15 \\
	12x_0-60x_1+65x_2 & = & \sqrt{12^2+60^2+65^2}\csc B-12 \\
	288x_0+105x_1-500x_2 & = & \sqrt{288^2+105^2+500^2}\csc C-288
\end{eqnarray*}
Their intersection is
\[
	\left({319\sqrt{133632161}\over1113945}-1,-{29\sqrt{133632161}\over1113945},{87\sqrt{133632161}\over1856575}\right).
\]
We want a central projection of this point on the sphere centered at $(-1,0,0)$, which is
\[
	\left({55\over\sqrt{3131}}-1,-{5\over\sqrt{3131}},{9\over\sqrt{3131}}\right).
\]
Finally we want to relocate this point to its stereographic projection, which is
\[
	\left(-{5\over53}\left(\sqrt{3131}-55\right),{9\over53}\left(\sqrt{3131}-55\right)\right).
\]
This point is shown at the left in Fig.~\ref{fig:centers5}.

\subsection{Converting homogeneous coordinates when $\mathbf{K=-1}$.}
As before, the projected vertices of the triangle of reference are $A(4/5,3/5)$, $B(-8/15,2/5)$, and $C(-2/5,-8/15)$, and the polar coordinates are $r_A=1$, $\theta_A=(4/5,3/5)$, $r_B=2/3$, $\theta_B=(-4/5,3/5)$, $r_C=2/3$, and $\theta_C=(-3/5,-4/5)$.  Then we can calculate the coordinates on the embedded sphere using our formulas.  We note that the first coordinate is imaginary when the sphere has negative curvature.  We have
\[
	A'=\left({2\over3}i,{16\over15},{4\over5}\right)
	\hspace{30pt}B'=\left({1\over4}i,-{3\over5},{9\over20}\right)
	\hspace{30pt}C'=\left({1\over4}i,-{9\over20},-{3\over5}\right).
\]
Since the center $O$ of the embedded sphere is at $(-\sqrt{1/K},0,0)=(-i,0,0)$, the three planes in which the triangle edges are embedded are
\begin{eqnarray*}
	OB'C': & & 3ix_0-7x_1-x_2 = 3 \\
	OC'A': & & 84ix_0+600x_1-625x_2 = 84 \\
	OA'B': & & 288ix_0-75x_1+700x_2 = 288
\end{eqnarray*}
The normals to these planes allow us to compute the cosines of the angles of the triangle of reference, which are
\[
	\cos A={506692\over\sqrt{306856798489}}
	\hspace{36pt}\cos B={1039\over\sqrt{16919921}}
	\hspace{36pt}\cos C={3827\over\sqrt{30486329}}
\]
From these we can compute the sines of the angles.
\[
	\sin A={223875\over\sqrt{306856798489}}
	\hspace{36pt}\sin B={3980\over\sqrt{16919921}}
	\hspace{36pt}\sin C={3980\over\sqrt{30486329}}
\]
We translate the planes above by the homogeneous coordinates, which we shall take as those for the intersection of the medians, ($\csc A:\csc B:\csc C$).  The new planes are
\begin{eqnarray*}
	3ix_0-7x_1-x_2 & = & \sqrt{-3^2+7^2+1^2}\csc A+3 \\
	84ix_0+600x_1-625x_2 & = & \sqrt{-84^2+600^2+625^2}\csc B+84 \\
	288ix_0-75x_1+700x_2 & = & \sqrt{-288^2+75^2+700^2}\csc C+288
\end{eqnarray*}
Their intersection is
\[
	\left(\left({\sqrt{12581128738049}\over712818}-1\right)i,{\sqrt{12581128738049}\over178204500},{13\sqrt{12581128738049}\over59401500}\right).
\]
The central projection of this point on the sphere centered at $(-i,0,0)$ is
\[
	\left(\left({250\over\sqrt{60978}}-1\right)i,{1\over\sqrt{60978}},{117\over\sqrt{60978}}\right).
\]
Finally we want to relocate this point to its stereographic projection, which is
\[
	\left({1\over761}\left(250-\sqrt{60978}\right),{39\over761}\left(250-\sqrt{60978}\right)\right).
\]
This point is shown at the right in Fig.~\ref{fig:centers5}.

\section{The median point}

\begin{figure}[t]
	\begin{center}
	\begin{picture}(400,200)(10,0)
		\includegraphics[width=6in]{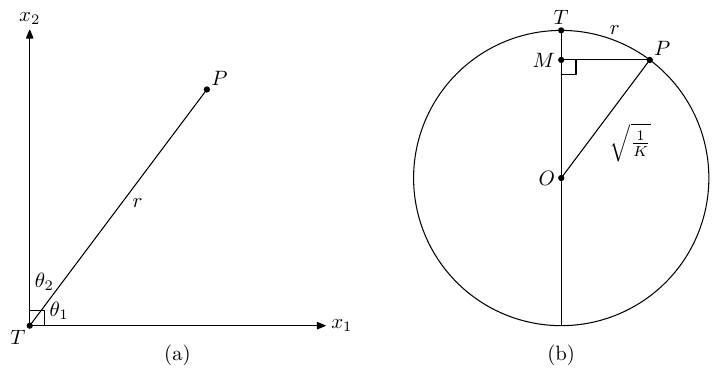}
	\end{picture}
	\end{center}
	\caption{The point $P$ on the left is given the coordinates $r\theta$, where $r$ is the distance of $P$ from an origin point $T$, and $\theta$ is an ordered pair of direction cosines, here $(\cos\theta_1,\cos\theta_2)$.  On the right, we see arc $TP$ in a cross section of the embedded sphere with Gaussian curvature $K$.  Note that $MP=\sing(r)$ and $OM=\sqrt{1/K}\cosg(r)$.}
	\label{fig:centers4}
\end{figure}

In this section we discuss an important non-Euclidean center point that we call the {\em median point}.  In our list of triangle centers below, we provide homogeneous coordinates for the median points of the vertices, edges, and interior of a triangle.  To define the median point, we first assign coordinates of the form $r\theta$ to points in a space of constant curvature $K$.  We select an origin in the space and let $r$ represent the distance from the origin of a point in the space.  We select a set of orthogonal rays from the origin and let $\theta$ represent the set of direction cosines determined by these rays.  As shown in Fig.~\ref{fig:centers4}, let the origin be $T$ and the rays be $x_1$ and $x_2$.  (There would be more rays if the space had more than two dimensions.)  If we draw a line $TP$ to any point $P$ in the space, the coordinates will be $r\theta$, where $r$ is the length of $TP$, and $\theta$ is the cosines of the angles between $TP$ and each of the rays.  As shown in Fig.~\ref{fig:centers4}(a), the coordinates of $P$ are $r(\cos\theta_1,\cos\theta_2)$.

Consider a set of $n$ points $P_i$, where $1\leq i \leq n$.  The coordinates of $P_i$ are $r_i\theta_i$.  If $0=\sum_i \sing(r_i)\theta_i$, we shall say that the coordinate origin point $T$ is the {\em median point} of that set of points.  If there is just one point, it must be its own median point.  If there are two points, they must be equidistant from $T$ in opposite directions, so that the median point of two points is the midpoint of the line joining them.  When there are more points, we have a simple process for determining their median point.

As shown in Fig.~\ref{fig:centers4}(b), if we embed our curved space in a Euclidean space, there is a simple way of transforming our $r\theta$ coordinates for our space of curvature $K$ into Euclidean coordinates of the space in which the sphere is embedded.  The cross section of the sphere includes the center of the sphere $O$, the coordinate origin point $T$, and a particular point $P$.  We see that
\[
	{MP\over OP}=\sin{r\sqrt{K}},
\]
and thus
\[
	MP={\sin{r\sqrt{K}}\over\sqrt{K}}=\sing r.
\]
Similarly, $OM=\sqrt{1/K}\cosg(r)$.

We let $O$ be the origin for our Euclidean coordinates with the direction from $O$ to $T$ being the $x_0$ axis and the other axes having the same orientation at $T$ as they do on the sphere.  Then the Euclidean coordinates for $r\theta$ on the sphere become $(\sqrt{1/K}\cosg(r),\sing(r)\theta)$.  For example, if $K=1$ and $r\theta={\pi\over 4}(\cos{\pi\over6},\cos{\pi\over3})$, our Euclidean coordinates would be
\[
	\left({\cosg r\over\sqrt{K}},\sing r\left(\cos{\pi\over6},\cos{\pi\over3}\right)\right)=\left({\sqrt{2}\over2},{\sqrt{6}\over4},{\sqrt{2}\over4}\right).
\]

Now we consider a set of points $P_i$ such that their median point is $T$, which by definition means that
\[
	0=\sum_i\sing(r_i)\theta_i.
\]
The Euclidean coordinates for these points become $(\sqrt{1/K}\cosg(r_i),\sing(r_i)\theta_i)$.  Let us find the Euclidean centroid of these $n$ points.  It is simply
\[
	\left({\sum_i\cosg r_i\over n\sqrt{K}},{1\over n}\sum_i\sing(r_i)\theta_i)\right)
		=\left({\sum_i\cosg r_i\over n\sqrt{K}},0\right).
\]
In other words, all but the first Euclidean coordinate is zero, so the Euclidean centroid must lie on the line $OT$.  Thus, if we have a set of points in a curved space that we embed in a Euclidean space, the median point of those points is found by finding the Euclidean centroid of the points (a point within the sphere), and then projecting it onto the sphere from the center of the sphere.  Of course this won't work if the Euclidean centroid is in fact the center of the sphere, but in that case we say that the median point is not defined.

\begin{figure}
	\begin{center}
	\begin{picture}(400,200)(0,0)
		\includegraphics[width=6in]{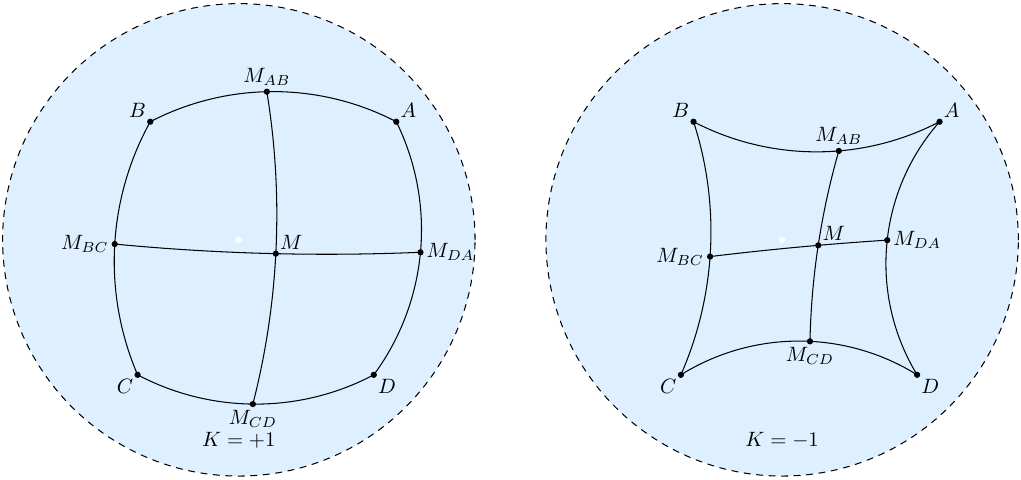}
	\end{picture}
	\end{center}
	\caption{To find the median point $M$ of points $A$, $B$, $C$, and $D$, we find the intersection of the median joining the midpoints of edges $AB$ and $CD$ with the median joining the midpoints of edges $BC$ and $DA$.}
	\label{fig:centers8}
\end{figure}

This use of the Euclidean centroid gives us a method for finding a median point.  We know that the centroid of a set of points must lie on a line joining the centroids of two nonempty subsets of these points when each point lies in one of the subsets.  Thus the centroid of the vertices of a Euclidean triangle must lie on a line joining one vertex to the midpoint of the edge joining the other two vertices.  But if the median point of a set of points lies on a central projection of the centroid, the same relationship holds for the median points of the subsets, since a central projection preserves straight lines.  We call this central projection of the line joining the median points of the two subsets a {\em median}.  Thus, the median point of a set of points in a space of constant curvature must lie on a median joining the median points of two complementary nonempty subsets of those points.  In other words, we can locate the median point of a non-Euclidean triangle at the intersection of its medians.  For four points, we take a median joining the midpoints of two pairs and repeat that with a different set of pairs as shown in Fig.~\ref{fig:centers8}.  This will also work to find the median point of the vertices of a non-Euclidean tetrahedron.

\section{A planar center of rotation}
If we consider a rigid set of points in a non-Euclidean plane of constant Gaussian curvature $K$, we can determine if a point in that plane is a center of rotation for these points as they rotate in the plane.  Using the same polar coordinate system that we used to define the median point, we say that the coordinate origin $T$ is a planar center of rotation of a rigid set of points if
\[
	0=\sum_i\sing(2r_i)\theta_i.
\]
Why is this point a center of rotation?  Lamphere~\cite{Lamphere} has shown that the centrifugal force of a particle rotating in non-Euclidean space of constant Gaussian curvature $K$ is equivalent to
\[
	{mv^2\over\tang r}.
\]
In non-Euclidean space the circumference of a circle of radius $r$ is $2\pi\sing r$, so that the velocity of objects rotating about a fixed point at a fixed number of revolutions per unit time is proportional to $\sing r$.  Assuming that our points have equivalent masses, the centrifugal force of a point is proportional to
\[
	{\sing^2r\over\tang r}={\sing r\cosg r}={\sing(2r)\over 2}.
\]
So the centrifugal force is proportional to $\sing(2r)$, and in order for these forces to cancel each other out at the coordinate origin $T$, the value of $\sum_i\sing(2r_i)\theta_i$ must be zero.

If $K$ is negative, this point is unique.  When $K$ is positive, however, this is not the case.  Consider a single point.  One center of rotation is the point itself.  Another center of rotation is any point on the polar of the point.  For two points a fixed distance apart, one center of rotation is the midpoint of the line connecting them.  Other centers of rotation include the poles of that line and the intersections of that line and the polar of the midpoint. 

Non-Euclidean rotations in more than two dimensions are complex.  In Euclidean space all axes of rotation of a rigid object pass through a common point.  In non-Euclidean space, this apparently is not true.  See Gunn~\cite{Gunn}.

\begin{figure}
	\begin{center}
	\begin{picture}(400,400)(14,0)
		\includegraphics[width=6in]{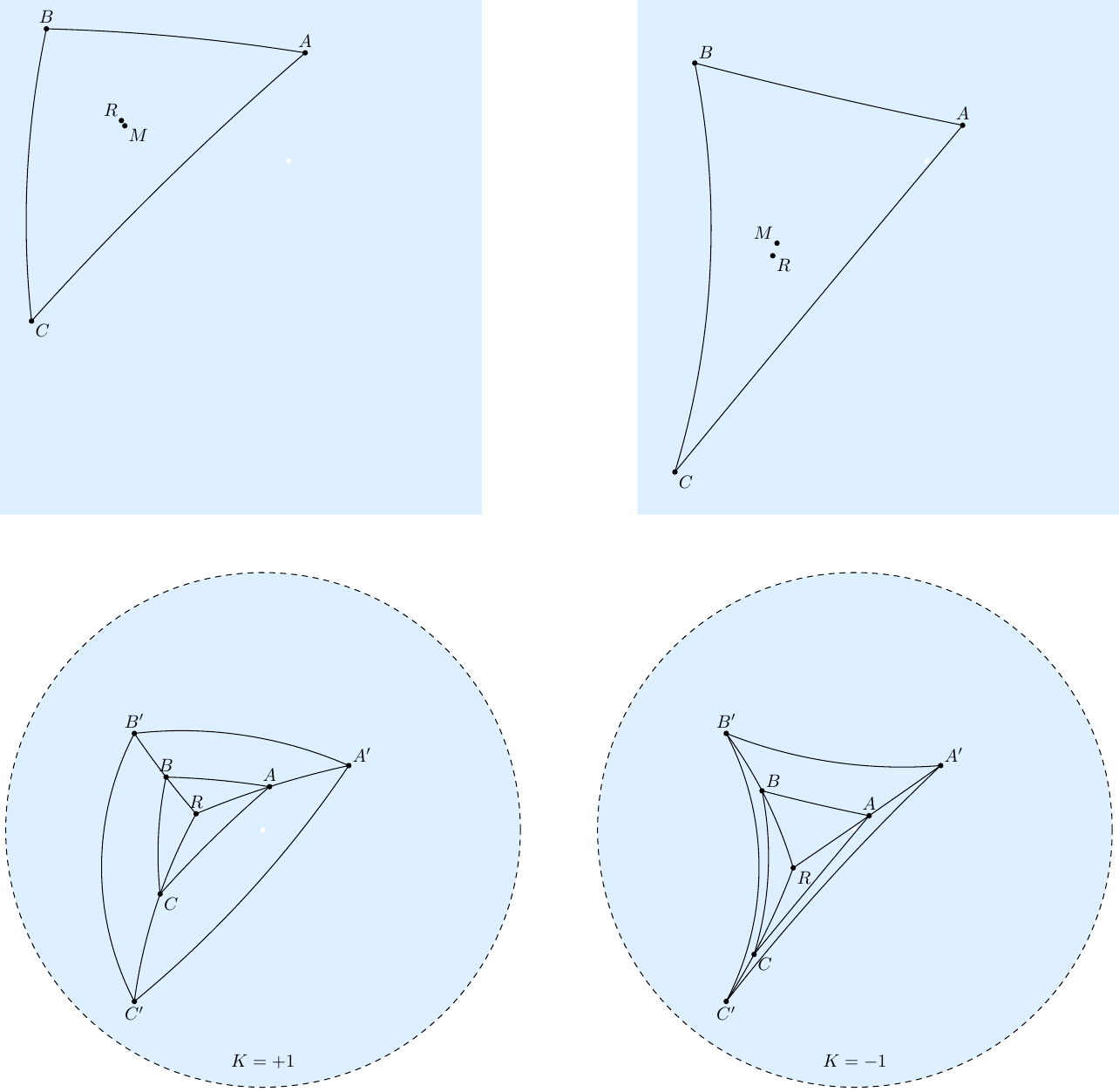}
	\end{picture}
	\end{center}
	\caption{Triangle $ABC$ is shown with the interior center of rotation $R$ of its vertices and the median point $M$ of its vertices.  If we extend $RA$ so that $RA'=2RA$ and similarly for the other vertices, $R$ is the median point of the vertices of triangle $A'B'C'$.}
	\label{fig:centers3}
\end{figure}

In Fig.~\ref{fig:centers3}, we show the median point and the interior center of rotation for triangles in spaces of positive and negative curvature.  As shown, if we extend the ray from the center of rotation to each point to double its length, the median point of the set of new points at the end of the extended rays will coincide with that interior center of rotation.

\begin{figure}
	\begin{center}
	\begin{picture}(400,200)(14,0)
		\includegraphics[width=6in]{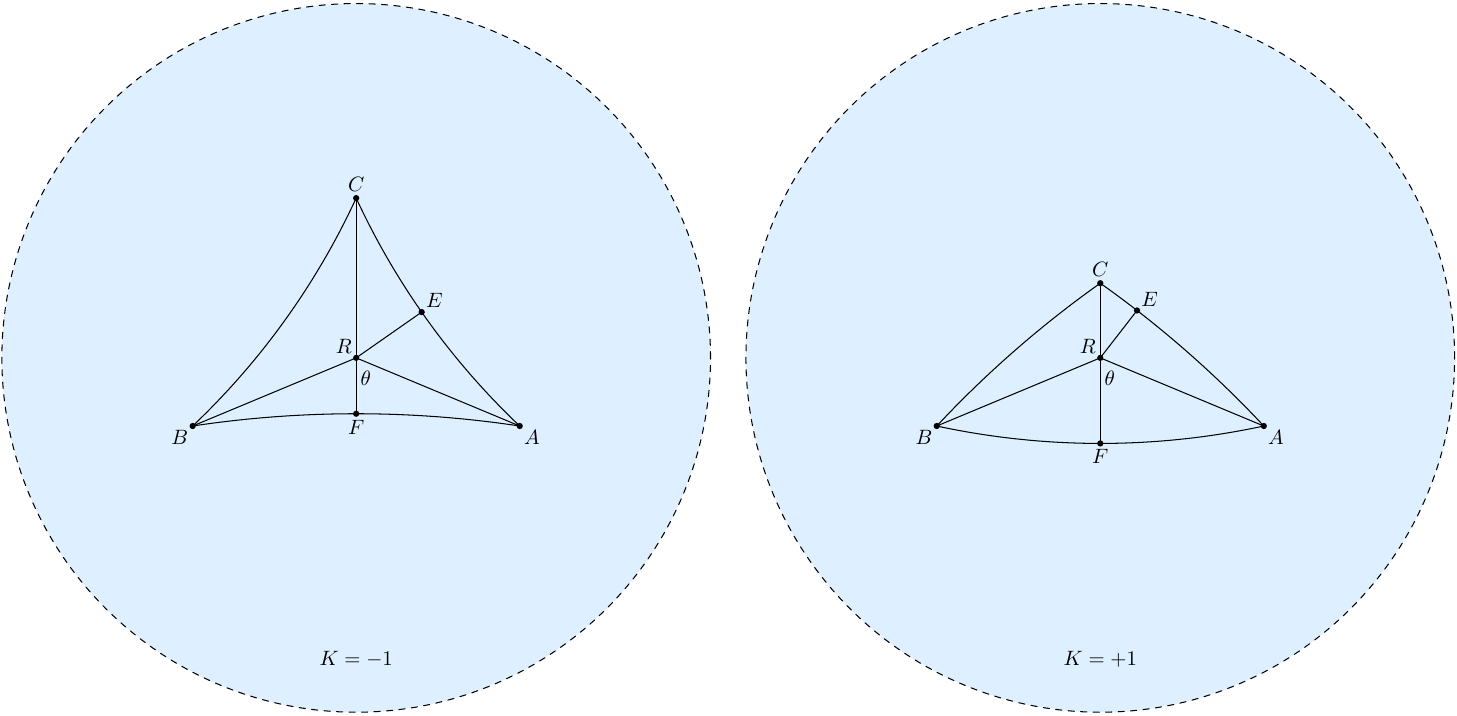}
	\end{picture}
	\end{center}
	\caption{In isosceles triangle $ABC$ with apex $C$, the interior center of rotation of its vertices is $R$.  The homogeneous coordinates of $R$ will be $(\sing RE : \sing RE : \sing RF)$.}
	\label{fig:centroid1}
\end{figure}

\subsection{Interior center of rotation of the vertices of an isosceles triangle}

Let us determine the homogeneous coordinates of the interior center of rotation of the vertices of an isosceles triangle $ABC$ with apex $C$ as shown in Fig.~\ref{fig:centroid1}.  Note that when the Gaussian curvature is positive, there are exterior centers of rotation at the poles of $CF$.  The homogeneous coordinates will be $(\sing RE : \sing RE : \sing RF)$.  First we determine $RF$.  We know that the three vectors must cancel out at point $R$, which requires that
\begin{equation}
	\sing 2RC=2\cos\theta\sing 2RA,
	\label{eqn:verts}
\end{equation}
where $\theta=\angle ARF$.  Note that $FA$ is half of $AB$ so that
\[
	2\cosg^2FA = 1+\cosg AB.
\]
In right triangle $RFA$,
\[
	\cos\theta={\tang RF\over\tang RA}
\]
and
\[
	\cosg RA = \cosg RF\cosg FA
\]
so that
\begin{eqnarray*}
	\sing 2RC & = & 4{\tang RF\over\tang RA}\sing RA\cosg RA \\
		& = & 4\tang RF\cosg^2RA \\
		& = & 4\tang RF\cosg^2 RF\cosg^2FA \\
		& = & \sing2RF(1+\cosg AB).
\end{eqnarray*}
From right triangle $CFA$, we have
\[
	\cosg CF = {\cosg AC\over\cosg FA},
\]
so that
\[
	\cosg 2CF = 2\cosg^2CF-1 = {4\cosg^2AC\over1+\cosg AB}-1.
\]
Also
\begin{eqnarray*}
	\cosg 2CF &=& \cosg2RC\cosg2RF-K\sing2RF\sing2RC \\
	{4\cosg^2AC\over1+\cosg AB}-1 & = & \cosg2RC\cosg2RF - K\sing^22RF(1+\cosg AB) \\
		& = & \cosg2RC\cosg2RF - (1-\cosg^22RF)(1+\cosg AB) \\
		& = & \sqrt{1-K\sing^22RC}\cosg2RF - (1-\cosg^22RF)(1+\cosg AB) \\
		& = & \sqrt{1-K\sing^22RF(1+\cosg AB)^2}\cosg2RF \\
		& & \hspace{20pt} - (1-\cosg^22RF)(1+\cosg AB) \\
		& = & \sqrt{1-(1-\cosg^22RF)(1+\cosg AB)^2}\cosg2RF \\
		& & \hspace{20pt} - (1-\cosg^22RF)(1+\cosg AB).
\end{eqnarray*}
The only unknown is $\cosg2RF$.  Solving, we obtain (the positive root)
\[
	\cosg2RF={\cosg AB+\cosg^2AB+4\cosg^2AC\over(1+\cosg AB)\sqrt{\cosg^2AB+8\cosg^2AC}}.
\]
This allows us to calculate either $\sing RF$ or $\cosg RF$ from the formula 
\[
	\cosg2RF=1-2K\sing^2RF=2\cosg^2RF-1.
\]

We now turn our attention to $RE$.  By the law of sines
\[
	\sing RE = {\sing{AB\over2}\sing RC\over\sing AC}={\sqrt{(1-\cosg AB)/(2K)}\sing RC\over\sing AC}.
\]
Also
\[
	\sing RC = \sing(CF-RF)=\cosg RF\sing CF-\sing RF\cosg CF.
\]
In right triangle $AFC$
\[
	\cosg CF = {\cosg AC\over\cosg FA}={\sqrt2\cosg AC\over\sqrt{1+\cosg AB}}.
\]
This allows us to calculate $\sing CF$ and therefore $\sing RE$ as functions of $AB$ and $AC$ (and $K$), which we leave to the reader.

We are now in a position to calculate the homogeneous coordinates of $R$ by evaluating $\sing RE/\sing RF$.  We have
\begin{eqnarray*}
	{\sing RE\over\sing RF} & = & {\sqrt{(1-\cosg AB)/(2K)}\left(\cosg RF\sing CF-\sing RF\cosg CF\right)
		\over\sing AC\sing RF} \\
	& = & {\sqrt{(1-\cosg AB)/(2K)}\left(\sing CF/\tang RF-\cosg CF\right)
		\over\sing AC} \\
	& = & {\sqrt{(1-\cosg AB)/(2K)}\left(\sqrt{(\cosg AB+1-2\cosg^2AC)/K}/\tang RF-\sqrt2\cosg AC\right)
		\over\sing AC\sqrt{1+\cosg AB}} \\
	& = & {\sing AB\left(\sqrt{(\cosg AB-\cosg 2AC)/(2K)}/\tang RF-\cosg AC\right)
		\over\sing AC(1+\cosg AB)}.
\end{eqnarray*}

Now let's work on $1/\tang RF$.  We have
\begin{eqnarray*}
	{1\over\tang RF} & = & {\cosg RF\over\sing RF}=\sqrt{K{1+\cosg2RF\over1-\cosg2RF}} \\
		& = & \sqrt{K{(1+\cosg AB)\left(\sqrt{8\cosg^2AC+\cosg^2AB}-\cosg AB\right)+4\cosg^2AC
			\over(1+\cosg AB)\left(\sqrt{8\cosg^2AC+\cosg^2AB}-\cosg AB\right)-4\cosg^2AC}} \\
		& = & {\sqrt{8K}\cosg AC\sqrt{\cosg AB-\cosg2AC}
		\over(1+\cosg AB)\left(\sqrt{8\cosg^2AC+\cosg^2AB}-\cosg AB\right)-4\cosg^2AC}.
\end{eqnarray*}

Plugging this into our earlier equation for $\sing RE/\sing RF$ and simplifying, we get
\[
	{\sing RE\over\sing RF}= {\sing AB\over\sing AC}\hspace{3pt}
		{\cosg AC\left(\sqrt{8\cosg^2AC+\cosg^2AB}-\cosg AB-2\right)
		\over4\cosg^2AC-(1+\cosg AB)\left(\sqrt{8\cosg^2AC+\cosg^2AB}-\cosg AB\right)}.
\]

Finally we have the homogeneous coodinates of $R$.  Substituting $b$ for $AC$ and $c$ for $AB$, the value of $(\sing RE:\sing RE:\sing RF)$ is
\[
	\left(\csc A\left(\cosg b\left(\sqrt{8\cosg^2b+\cosg^2c}-\cosg c-2\right)\right)\right.
\]
\[
	:\csc B\left(\cosg b\left(\sqrt{8\cosg^2b+\cosg^2c}-\cosg c-2\right)\right)
\]
\[
	:\left.\csc C\left(4\cosg^2b-(1+\cosg c)\left(\sqrt{8\cosg^2b+\cosg^2c}-\cosg c\right)\right)\right).
\]
We can convert this to functions of only the angles of the triangle with a law of cosines.
\[
	\cosg b = {\cos B+\cos A\cos C\over\sin A\sin C}\hspace{30pt}
	\cosg c = {\cos C+\cos A\cos B\over\sin A\sin B}.
\]

Having found the homogeneous coordinates of $R$, the interior planar center of rotation for the vertices of an isosceles triangle, we are still far from obtaining the homogeneous coordinates for the same center of a scalene triangle, which we pose as (presumably open) problem 4 below.

\subsection{Interior center of rotation of the edges of an isosceles triangle}

\begin{figure}
	\begin{center}
	\begin{picture}(400,200)(14,0)
		\includegraphics[width=6in]{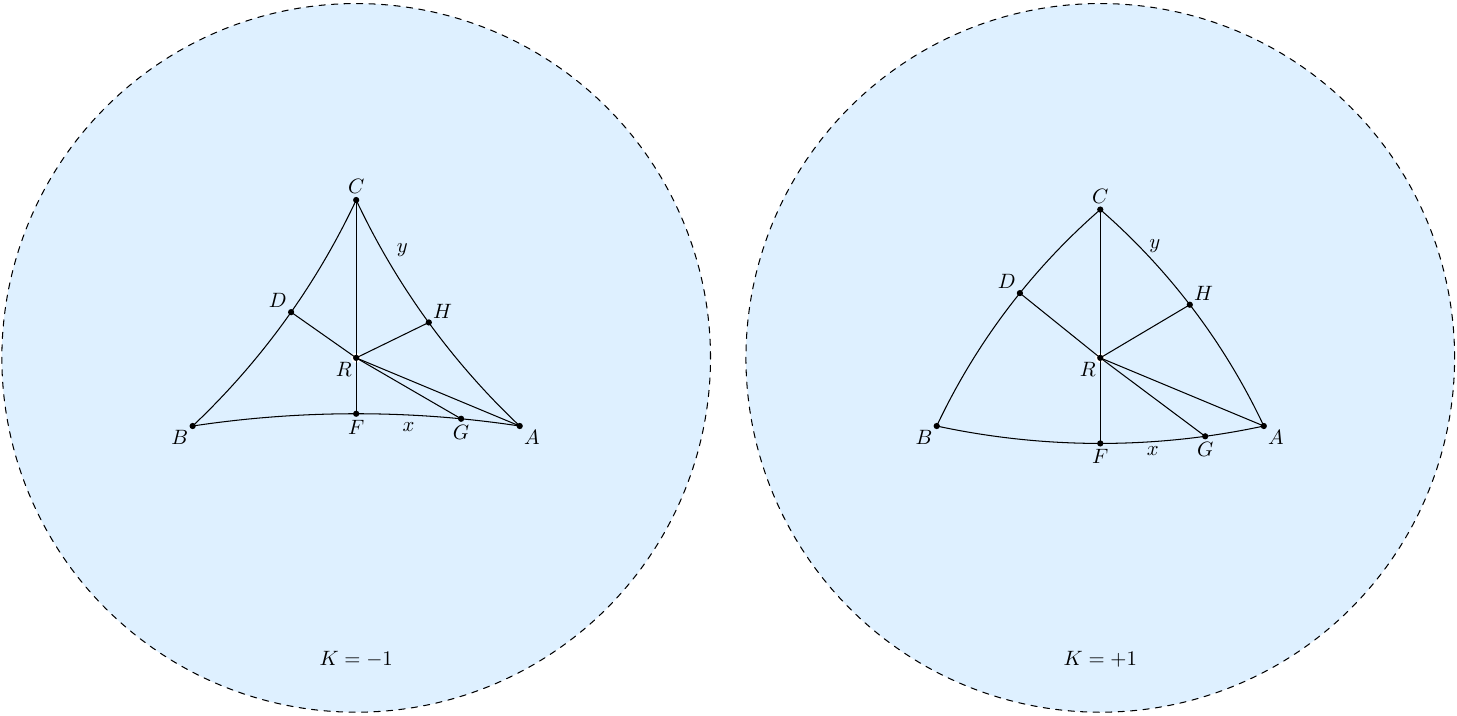}
	\end{picture}
	\end{center}
	\caption{In isosceles triangle $ABC$ with apex $C$, the interior center of rotation of its edges is $R$.  The homogeneous coordinates of $R$ will be $(\sing RD : \sing RD : \sing RF)$.}
	\label{fig:centroid2}
\end{figure}

We want to determine the rotation center of the edges of a non-Euclidean isosceles triangle.  Let the edges meeting at vertex $C$ of triangle $ABC$ be equal and $F$ be the foot of the altitude from $C$.  There is a point $R$ on $CF$ that is the interior center of rotation of the edges of the triangle.  As shown in Fig.~\ref{fig:centroid2}, let $G$ be a point on $FA$ and distance $FG$ be $x$, and let $H$ be a point on $AC$, and distance $CH$ be $y$.  Then $R$ is the center of rotation of the edges of $ABC$ if
\begin{equation}
	\int_0^{c/2}\sing2RG\cos FRG \;dx=\int_0^b\sing 2RH\cos CRH \;dy.
	\label{eqn:integral}
\end{equation}
This is because the component of $RG$ in the direction of $RF$ must be balanced by the component of $RH$ in the direction of $RC$.

We evaluate the left side of (\ref{eqn:integral}) first.  From right triangle $RFG$, we have
\[
	\cosg RG=\cosg RF \cosg x,
\]
\[
	\sin FRG={\sing x\over\sing RG}
\]
and
\[
	\tan FRG = {\tang x\over\sing RF}.
\]
Then
\begin{eqnarray*}
	\sing2RG & = & 2\sing RG\cosg RG= 2\cosg RF\cosg x\sing RG \\
		& = & {2\cosg RF\cosg x\sing x\over\sin FRG} = {\cosg RF\sing 2x\over\sin FRG}
\end{eqnarray*}
Adding in the other factor on the left-hand side of (\ref{eqn:integral}),
\begin{eqnarray*}
	\sing2RG\cos FRG & = & {\cosg RF\sing 2x\cos FRG\over\sin FRG} =  {\cosg RF\sing 2x\over\tan FRG} \\
		& = & {\sing RF\cosg RF\sing 2x\over\tang x} = {\sing2RF\cosg^2x}.
\end{eqnarray*}
Now we can integrate the left-hand side of (\ref{eqn:integral}) to obtain
\[
	\int_0^{c/2}\sing2RG\cos FRG \;dx = \sing 2RF\int_0^{c/2}\cosg^2x \;dx = \sing 2RF(c+\sing c)/4.
\]
We want to replace $\sing2RF$ with a function of $RC$.  Since $2RF=2CF-2RC$ and $\sing CF=\sin A\sing b$ and $\cosg CF = \cosg b/\cosg c/2$, we have
\begin{eqnarray*}
	\sing 2RF & = & \sing2CF\cosg2RC-\cosg2CF\sing2RC \\
		& = & 2\sing CF\cosg CF\cosg2RC-(2\cosg^2CF-1)\sing2RC \\
		& = & {\sin A\sing2b\over\cosg c/2}\cosg2RC
			-\left({2\cosg^2b\over\cosg^2{c/2}}-1\right)\sing2RC. \\
\end{eqnarray*}

Now let us deal with the integral on the right-hand side of (\ref{eqn:integral}).  We no longer have the luxury of a right triangle.  But we do know that
\[
	\cosg RH=\cosg RC\cosg y+K\sing RC\sing y\cos{C\over2}
\]
and
\begin{eqnarray*}
	\sing 2RH\cos CRH & = & 2\sing RH\cosg RH\cos CRH.
\end{eqnarray*}

We want to find $\sing RH\cos CRH$ in terms of $y$ and values independent of $y$.  Using a law of cosines,
\begin{eqnarray*}
	\sing RH\cos CRH & = & {\cosg y-\cosg RC\cosg RH\over K\sing RC} \\
		& = & {(1-\cosg^2RC)\cosg y-\cosg RC\cosg RH+\cosg^2RC\cosg y\over K\sing RC} \\
		& = & {K\sing^2RC\cosg y-\cosg RC(\cosg RH-\cosg RC\cosg y)\over K\sing RC} \\
		& = &  \sing RC\cosg y-\cosg RC\sing y{\cosg RH-\cosg RC\cosg y\over K\sing RC\sing y} \\
		& = & \sing RC\cosg y-\cosg RC\sing y\cos{C\over2}. \\
\end{eqnarray*}
Then
\begin{eqnarray*}
	\sing2RH\cos CRH & = & 2\sing RH\cosg RH\cos CRH \\
		& = & 2\cosg RH\left(\sing RC\cosg y-\cosg RC\sing y\cos{C\over2}\right) \\
		& = & 2\left(\cosg RC\cosg y+K\sing RC\sing y\cos{C\over2}\right) \\
		& & \left(\sing RC\cosg y-\cosg RC\sing y\cos{C\over2}\right).
\end{eqnarray*}

Since $RC$ and $C$ are independent of $y$, we can integrate to obtain
\begin{equation*}
\begin{split}
	& \int_0^b\sing 2RH\cos CRH\;dy = \\
	& {1\over8} \left(4b\sin^2{C\over2}+(3+\cos C)\sing2b\right) \sing2RC
		-\cos{C\over2}\sing^2b\cosg2RC.
\end{split}
\end{equation*}

Now we can replace both integrals in (\ref{eqn:integral}) to obtain
\begin{equation*}
\begin{split}
	& {1\over8} \left(4b\sin^2{C\over2}+(3+\cos C)\sing2b\right) \sing2RC
		-\cos{C\over2}\sing^2b\cosg2RC \\
	& = \sing 2RF(c+\sing c)/4 \\
	& = \left({\sin A\sing2b\over\cosg c/2}\cosg2RC
			-\left({2\cosg^2b\over\cosg^2{c/2}}-1\right)\sing2RC\right)(c+\sing c)/4.
\end{split}
\end{equation*}
This can be resolved into
\begin{eqnarray*}
	\tang2RC & = & {{2\sin A\sing2b(c+\sing c)\over\cosg c/2}+8\cos{C\over2}\sing^2b\over
	 	\left(4b\sin^2{C\over2}+(3+\cos C)\sing2b\right)
	 	+2\left({2\cosg^2b\over\cosg^2{c/2}}-1\right)(c+\sing c)} \\
		& = & {4\cos{C\over2}\sing^2b+\sing2b(c+\sing c)\sqrt{1-{\tang^2{c\over2}\over\tang^2b}}/\cosg{c\over2}
		\over(1-2\cosg2b+\cosg c)/\tang b+2b{\sing^2{c\over2}\over\sing^2b}
		+\left(2{\cosg^2b\over cosg^2{c\over2}}-1\right)(c+\sing c)}.
\end{eqnarray*}
This allows us to determine $\tang RC$ as
\[
	\tang RC = {\sqrt{1+K\tang^22RC}-1\over K\tang2RC}.
\]

Once we have the value of $\tang RC$, we can calculate $\sing RF$ and $\sing RD$.  These are
\begin{eqnarray*}
	\sing RF & = & \sing(CF-RC) = \sing CF\cosg RC - \cosg CF\sing RC \\
		& = & \sing b\sing B\cosg RC - \cosg b\sing RC/\cosg{c\over2} \\
		& = & {\sing b\sing B- \cosg b\tang RC/\cosg{c\over2}\over\cosg RC} \\
		& = & \sqrt{1+K\tang^2RC}\left(\sing b\sing B- \cosg b\tang RC/\cosg{c\over2}\right). \\
	\sing RD  & = & \sin{C\over2}\sing RC = {\sin{C\over2}\tang RC\over\sqrt{1+K\tang^2RC}}.
\end{eqnarray*}
The ratio is then
\[
	{\sing RD\over\sing RF} = {\sin{C\over2}\tang RC\over\left(\sin B\sing b-\cosg b\tang RC/\cosg{c\over2}\right)\left(1+K\tang^2RC\right)},
\]
and the homogeneous coordinates of the center of rotation are
\[
	\left(\sin{C\over2}\tang RC : \sin{C\over2}\tang RC : \left(\sin B\sing b-\cosg b\tang RC/\cosg{c\over2}\right)\left(1+K\tang^2RC\right)\right).
\]

\subsection{Interior center of rotation of the interior of an isosceles triangle}

\begin{figure}
	\begin{center}
	\begin{picture}(400,200)(14,0)
		\includegraphics[width=6in]{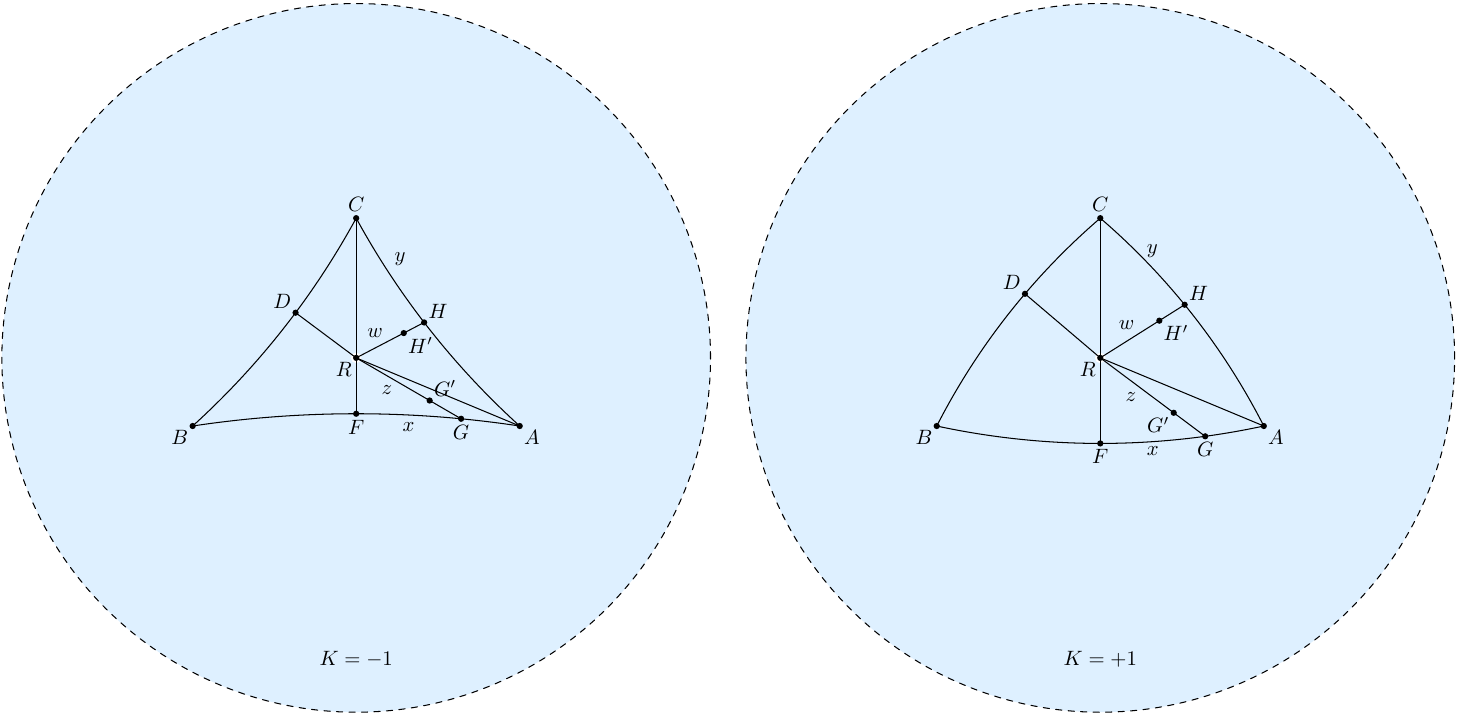}
	\end{picture}
	\end{center}
	\caption{In isosceles triangle $ABC$ with apex $C$, the center of rotation of its interior is $R$.  The homogeneous coordinates of $R$ will be $(\sing RD : \sing RD : \sing RF)$.}
	\label{fig:centroid3}
\end{figure}

We want to determine the rotation center of the interior of a non-Euclidean isosceles triangle.  Let the edges meeting at vertex $C$ of triangle $ABC$ be equal and $F$ be the foot of the altitude from $C$.  There is a point $R$ on $CF$ that is the interior center of rotation of the interior of the triangle.  As shown in Fig.~\ref{fig:centroid3}, let $G$ be a point on $FA$ and distance $FG$ be $x$, and let $H$ be a point on $AC$, and distance $CH$ be $y$.  Since we shall need double integration to deal with the interior points of the triangle, we shall let $G'$ be a point on $RG$ with the distance $RG'$ being $z$.  Also $H'$ is a point on $RH$ with the distance $RH'$ being $w$.

We shall first integrate sectors centered at $R$ with bases along $FA$ and $CA$.  We then integrate the components of these sectors in the directions of $RC$ and $RF$ respectively.  These components should be equal when $R$ is the center of rotation.  We begin with the integration of the sectors.

\begin{figure}
	\begin{center}
	\begin{picture}(400,220)(0,0)
		\includegraphics[width=5.5in]{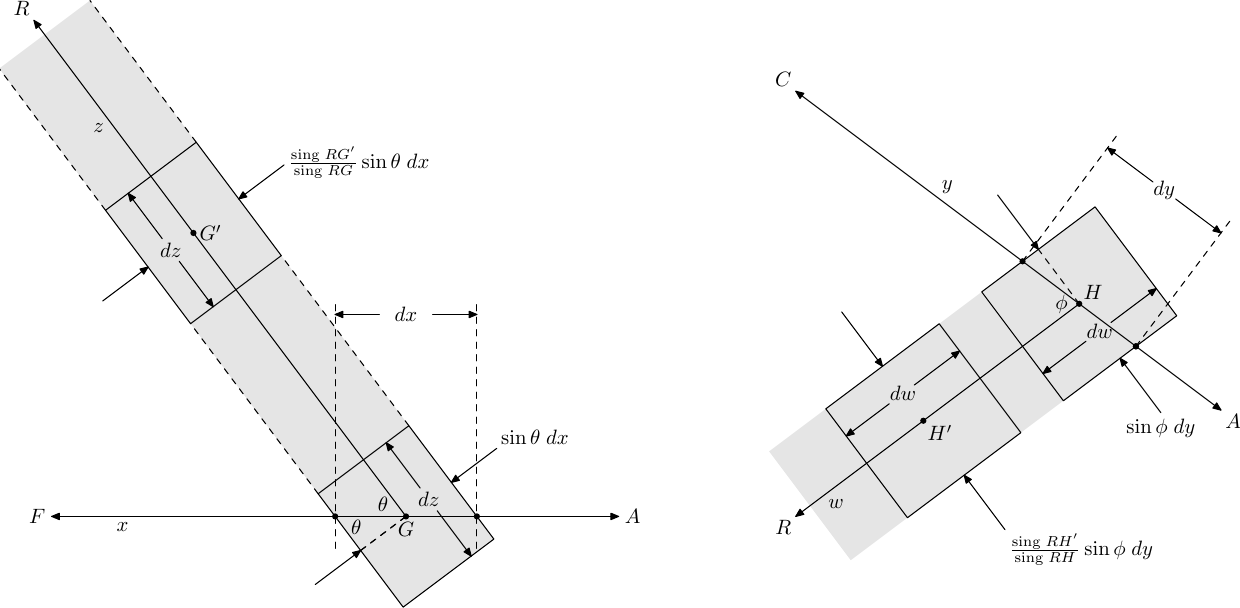}
	\end{picture}
	\end{center}
	\caption{The gray areas are parts of sectors for which we want to integrate the small rectangles around points $G'$ and $H'$ to determine the centrifugal force of the sector.}
	\label{fig:double_int}
\end{figure}

In Fig.~\ref{fig:double_int}, the gray area represents a portion of a sector of a circle centered at $R$ with radius $RG+dz/2$.  The intersection of the sector with edge $AB$ of the isosceles triangle is a segment of length $dx$ centered at $G$.  The width of the sector at $G$ is $\sin\theta\;dx$, where $\theta=\angle RGF$.  The width of the sector at $G'$ is ${\sing RG'\over\sing RG}\sin\theta\;dx$.  The centrifugal force $F_{dz}$ of the small rectangle around $G'$ is proportional to
\begin{eqnarray*}
	F_{dz} & = & \sing2RG'{\sing RG'\over\sing RG}\sin\theta\;dx\;dz
		=  \sing2z{\sing z\over\sing RG}\sin\theta\;dx\;dz \\
		& = & 2\sing z\cosg z{\sing z\over\sing RG}\sin\theta\;dx\;dz
		= 2{\sing^2z\cosg z\over\sing RG}\sin\theta\;dx\;dz \\
		& = & (\cosg z-\cosg^3z){2\sin\theta\over K\sing RG}\;dx\;dz.
\end{eqnarray*}

The centrifugal force $F_s$ of the entire sector is then the integral
\begin{eqnarray*}
	F_s & = & \int_0^{RG}(\cosg z-\cosg^3z){2\sin\theta\over K\sing RG}\;dx\;dz \\
		& = & {2\sin\theta\over K\sing RG}\;dx\int_0^{RG}(\cosg z-\cosg^3z)\;dz
		= {2\over3}\sing^2RG\sin\theta\;dx.
\end{eqnarray*}

The sector through $RH$ can similarly be determined to have a centrifugal force proportional to
\[
	{2\over3}\sing^2RH\sin\phi\;dy,
\]
where $\phi=\angle RHC$.

Now we can say that $R$ is the planar center of rotation of the points within the isosceles triangle $ABC$ if
\begin{equation}
	\int_0^{c/2}\sing^2RG\sin RGF\cos FRG \;dx=\int_0^b\sing^2RH\sin RHC\cos CRH \;dy.
	\label{eqn:doubleint}
\end{equation}
Let's first integrate the left-hand side of this.  We have
\[
	\sing^2RG\sin RGF\cos FRG = \sing^2RG{\sing RF\over\sing RG}{\tang RF\over\tang RG}
\]
\[
	= \cosg RG\sing RF\tang RF= \cosg x\cosg RF\sing RF\tang RF = \sing^2RF\cosg x.
\]
Integrating, we have
\[
	\int_0^{c/2}\sing^2RG\sin RGF\cos FRG \;dx = \sing^2RF\int_0^{c/2}\cosg x \;dx = \sing^2RF\sing{c\over2}.
\]

Now let's look at the right-hand side of (\ref{eqn:doubleint}).  We have
\[
	\sing^2RH\sin RHC\cos CRH = \sing^2RH{\sing RC\over\sing RH}\sin{C\over2}\cos CRH
\]
\begin{eqnarray*}
		& = & \sin{C\over2}\sing RC\sing RH\cos CRH \\
		& = & \sin{C\over2}{\sing RC\over2\cosg RH}\sing2RH\cos CRH \\
		& = & \sin{C\over2}{\sing RC\over2\cosg RH}2\cosg RH
			\left(\sing RC\cosg y-\cosg RC\sing  y\cos{C\over2}\right) \\
		& = & \sin{C\over2}\sing RC\left(\sing RC\cosg y-\cosg RC\sing  y\cos{C\over2}\right). \\
\end{eqnarray*}
Integrating, we have
\[
	\int_0^b\sing^2RH\sin RHC\cos CRH \;dy
\]
\begin{eqnarray*}
	& = & \sin{C\over2}\sing RC\int_0^b\sing RC\cosg y-\cosg RC\sing  y\cos{C\over2}\;dy \\
	& = & \sin{C\over2}\sing RC\left(\sing RC\sing b+{1\over K}\cosg RC(\cosg  b-1)\cos{C\over2}\right) \\
	& = & \sin{C\over2}\sing RC\left(\sing RC\sing b-{\cosg RC\sing ^2b\over1+\cosg b}\cos{C\over2}\right).
\end{eqnarray*}

Now we must replace $RF$ with $RC$ in the integral of the left-hand side of (\ref{eqn:doubleint}) in order to have the single unknown $RC$ on both sides.  Since $RF=CF-RC$, we have
\[
	\sing RF = \sing CF\cosg RC-\cosg CF\sing RC.
\]
Also $\sing CF = \sin B\sing b$ and $\cosg CF=\cosg b/\cosg{c\over2}$, so that
\[
	\sing RF = \sin B\cosg RC\sing b-\cosg b\sing RC/\cosg{c\over2}.
\]

We need to solve
\begin{multline*}
	\left(\sin B\cosg RC\sing b-\cosg b\sing RC/\cosg{c\over2}\right)^2\sing{c\over2} \\
	=\sin{C\over2}\sing RC\left(\sing RC\sing b-{\cosg RC\sing ^2b\over1+\cosg b}\cos{C\over2}\right).
\end{multline*}
Letting $\sing{c\over2}=\sing{b}\sin{C\over2}$, we divide both sides by $\sing{b}\sin{C\over2}\cosg^2RC$ to get
\[
	\left(\sin B\sing b-\cosg b\tang RC/\cosg{c\over2}\right)^2
	=\tang RC\left(\tang RC-{\sing{b}\over1+\cosg b}\cos{C\over2}\right),
\]
a quadratic equation with unknown $\tang RC=x$.
\begin{eqnarray*}
	0 & = & \left(1-{\cosg^2b\over\cosg^2{c\over2}}\right)x^2+\left({2\sin B\sing b\cosg b\over\cosg{c\over2}}
		-{\sing b\cos{C\over2}\over1+\cosg b}\right)x-\sin^2B\sing^2b \\
		& = & \left(1-\cosg^2CF\right)x^2+\left({2\sin B\sing b\cosg b\over\cosg{c\over2}}
		-{\sing b\cos{C\over2}\over1+\cosg b}\right)x-\sin^2B\sing^2b \\
		& = & K\sing^2CFx^2+\left({2\sin B\sing b\cosg b\over\cosg{c\over2}}
		-{\sing b\cos{C\over2}\over1+\cosg b}\right)x-\sin^2B\sing^2b \\
		& = & K\sin^2B\sing^2b\;x^2+\left({2\sin B\sing b\cosg b\over\cosg{c\over2}}
		-{\sing b\cos{C\over2}\over1+\cosg b}\right)x-\sin^2B\sing^2b \\
		& = & K\sin^2B\sing b\;x^2+\left({2\sin B\cosg b\over\cosg{c\over2}}
		-{\cos{C\over2}\over1+\cosg b}\right)x-\sin^2B\sing b. \\
\end{eqnarray*}

The solution is
\[
	\tang RC = {{\cos{C\over2}\over1+\cosg b}-{2\cosg b\sin B\over\cosg{c\over2}}
	+\sqrt{4\sin^2B+{\cosg^2{C\over2}\over\left(1+\cosg b\right)^2}
	-{4\cos{C\over2}\sin B\cosg b\over\cosg{c\over2}\left(1+\cosg b\right)}}
	\over2K\sin^2B\sing b}.
\]

Once we have the value of $\tang RC$, we can calculate $\sing RF$ and $\sing RD$.  These are
\begin{eqnarray*}
	\sing RF & = & \sing(CF-RC) = \sing CF\cosg RC - \cosg CF\sing RC \\
		& = & \sing b\sing B\cosg RC - \cosg b\sing RC/\cosg{c\over2} \\
		& = & {\sing b\sing B- \cosg b\tang RC/\cosg{c\over2}\over\cosg RC} \\
		& = & \sqrt{1+K\tang^2RC}\left(\sing b\sing B- \cosg b\tang RC/\cosg{c\over2}\right). \\
	\sing RD  & = & \sin{C\over2}\sing RC = {\sin{C\over2}\tang RC\over\sqrt{1+K\tang^2RC}}.
\end{eqnarray*}
The ratio is then
\[
	{\sing RD\over\sing RF} = {\sin{C\over2}\tang RC\over\left(\sin B\sing b-\cosg b\tang RC/\cosg{c\over2}\right)\left(1+K\tang^2RC\right)},
\]
and the homogeneous coordinates of the center of rotation are
\[
	\left(\sin{C\over2}\tang RC : \sin{C\over2}\tang RC : \left(\sin B\sing b-\cosg b\tang RC/\cosg{c\over2}\right)\left(1+K\tang^2RC\right)\right).
\]

\section{A short list of non-Euclidean triangle centers}
Here $S=(A+B+C)/2$ and $s=(a+b+c)/2$.  Note that many of the homogeneous coordinates are equivalent to the trilinear coordinates in Euclidean space.  The homogeneous coordinates of the Euler line~\cite{Akopyan} are
\[
	[\cos(2A\!-\!S)\sin(B\!-\!C):\cos(2B\!-\!S)\sin(C\!-\!A):\cos(2C\!-\!S)\sin(A\!-\!B)].
\]

\medskip

\noindent1.  The incenter, the center of the incircle.  $(1:1:1)$

\noindent2.  The vertex median point, the intersection of the medians.  $(\csc A:\csc B:\csc C)$

\noindent3.  The circumcenter, the center of the circumcircle.  $(\sin(S-A):\sin(S-B):\sin(S-C))$

\noindent4.  The orthocenter, the intersection of the altitudes.  $(\sec A:\sec B:\sec C)$

\noindent5.  The Euler circle center, the center of the circle that is externally tangent to the three excircles of the triangle.  $(\cos(B-C):\cos(C-A):\cos(A-B))$

\noindent6.  The symmedian point, the intersection of the symmedians.  $(\sin A:\sin B:\sin C)$

\noindent7.  The Gergonne point, the intersection of the cevians through the touch points of the incircle.
\[
	\left({\csc A\over\sing(s-a)}:{\csc B\over\sing(s-b)}:{\csc C\over\sing(s-c)}\right)
\]

\noindent8.  The Nagel point, the intersection of the cevians through the touch points of the excircles.
\[
	\left({\sing(s-a)\over\sin A}:{\sing(s-b)\over\sin B}:{\sing(s-c)\over\sin C}\right)
\]

\noindent9.  The mittenpunkt, the intersection of the lines through the midpoint of each edge and the center of the excircle on the other side of that edge from its opposite vertex.
\[
	(-\sin A+\sin B+\sin C:\sin A-\sin B+\sin C:\sin A+\sin B-\sin C)
\]

\noindent10.  The edge median point.
\[
	\left({\tang b+\tang c\over\sin A}:{\tang c+\tang a\over\sin B}:{\tang a+\tang b\over\sin C}\right)
\]

\noindent11.  The Feuerbach point, the common point of the incircle and the Euler circle (the circle that is externally tangent to the three excircles).
\[
	\left(1-\cos(B-C):1-\cos(C-A):1-\cos(A-B)\right)
\]

\noindent12.  The incenter of the medial triangle, with edges $a_m$, $b_m$, and $c_m$.  We can use the generalized law of cosines ($\cosg c=\cosg a\cosg b+K\sing a\sing b\cos C$) to evaluate these edges, replacing, e.g.,  $a$, $b$, and $c$ with $a/2$, $b/2$, and $c_m$.
\[
	\left({\cosg{b\over2}\sing c_m+\cosg{c\over2}\sing b_m\over\sing{a\over2}}:{\cosg{c\over2}\sing a_m+\cosg{a\over2}\sing c_m\over\sing{b\over2}}:{\cosg{a\over2}\sing b_m+\cosg{b\over2}\sing a_m\over\sing{c\over2}}\right)
\]

\noindent13.  The triangle median point, the median point of the points in the interior of the triangle.  We need to take the limit, $(\csc A:\csc B:\csc C)$, when $K=0$.
\[
	\left(a-b\cos C-c\cos B:b-c\cos A-a\cos C,c-a\cos B-b\cos A\right)
\]

\noindent14.  The polar median point, the median point of the points in the interior of the polar triangle, or the limit of this point when $K=0$.
\[
	\left({\pi-A\over\sin A}:{\pi-B\over\sin B}:{\pi-C\over\sin C}\right)
\]
\noindent15.  The cevian bisector center, the intersection of the three cevians that bisect the area of the triangle.
\[
	\left({1\over{\sin (S\!-\!A)-\cos A}}:{1\over{\sin (S\!-\!B)-\cos B}}:{1\over{\sin (S\!-\!C)}-\cos C}\right)
\]

\noindent16.  The pseudoaltitude center~\cite{Akopyan}, where the cevian $AD$ is a pseudoaltitude if 
\[
	\angle BDA - \angle ABD -\angle DAB = \angle ADC - \angle CAD-\angle DCA.
\]
The  homogeneous coordinates of the pseudoaltitude center are
\[
	\left({1\over{\sin (S\!-\!A)+\cos A}}:{1\over{\sin (S\!-\!B)+\cos B}}:{1\over{\sin (S\!-\!C)}+\cos C}\right).
\]

\noindent17.  The trisection center D such that  the areas of triangles $ABD$, $BCD$, and $CAD$ are equal.
\[
	\left(\csc(A-(2S\!-\!\pi)/3):\csc(B-(2S\!-\!\pi)/3):\csc(C-(2S\!-\!\pi)/3)\right).
\]

\section{Problems}
Here are some problems regarding the homogeneous coordinates of non-Euclidean triangle centers.  Perhaps only the first two have known solutions.

\noindent1.  What is the polar median point (triangle center 14 above) called in Euclidean geometry?

\noindent2.  Show that the limit when $K\rightarrow0$ of the generalized law of cosines (triangle center 12 above) is equivalent to $c^2=a^2+b^2-2ab\cos C$.

\noindent3.  What is the geometric description of the non-Euclidean triangle center with homogeneous coordinates $(\cos A:\cos B:\cos C)$? 

\noindent4.  What are the homogeneous coordinates of the interior planar center of rotation of the vertices of a non-Euclidean triangle?

\noindent5.  What are the homogeneous coordinates of the interior planar center of rotation of the edges of a non-Euclidean triangle?

\noindent6.  What are the homogeneous coordinates of the interior planar center of rotation of the interior of a non-Euclidean triangle?

\section{Barycentric coordinates for non-Euclidean triangles}

While Coxeter's homogeneous coordinates discussed here are trilinear coordinates in the Euclidean case, other authors use homogeneous coordinates that are barycentric coordinates in the Euclidean case.  G.Horv\`ath \cite{Ghorvath}  uses what is equivalent to the generalized polar sines \cite{Kokkendorff} of the triangles $BCM$, $CAM$, and $ABM$ as the coordinates for point $M$ in triangle $ABC$.  These coordinates can be obtained from Coxeter's $(x_0:x_1:x_2)$ by multiplying by the generalized sines of the edges to get $(x_0\sing a:x_1\sing b:x_2\sing c)$.


\end{document}